\documentclass{commtr}
\usepackage{appendix}
\usepackage{titlesec}
\usepackage{natbib}
\usepackage{algorithm}
\usepackage{algorithmic}
\usepackage{amsthm}
\newtheorem{definition}{Definition}
\newtheorem{assumption}{Assumption}
\newtheorem{theorem}{Theorem}
\newtheorem{corollary}{Corollary}
\usepackage{amsmath}
\usepackage{mathrsfs}
\usepackage{ragged2e}
\usepackage{comment}

\title{Model-Targeted Data Poisoning Attacks against ITS Applications with Provable Convergence}

\DeclareUnicodeCharacter{FF0C}{,}

\setcitestyle{aysep={,},yysep={,},notesep={;}}
\commtrsetup{
title    = {Model-Targeted Data Poisoning Attacks against ITS Applications with Provable Convergence },
author   = {Xin Wang$^{\rm a}$, Feilong Wang$^{\rm b}$, Yuan Hong $^{\rm c}$,  R. Tyrrell Rockafellar $^{\rm a}$, Xuegang (Jeff) Ban$^{\rm a,*}$},
address  = {
\tdd{a} Department of Civil and Environmental Engineering, University of Washington, Seattle, WA, 98195, United States \sep
\tdd{b} School of Transportation and Logistics, Southwest Jiaotong University, Chengdu, 610032, China \mbox{     }  (work conducted while a Ph.D. student at the University of Washington)\sep
\tdd{c}  School of Computing,  University of Connecticut, Storrs, CT, 06269, United States\sep
\tdd{*} Corresponding author.
},
email    = {banx@uw.edu},
abstract = {
The growing reliance of intelligent systems on data makes the systems vulnerable to data poisoning attacks. Such attacks could compromise machine learning or deep learning models by disrupting the input data. Previous studies on data poisoning attacks are subject to specific assumptions, and limited attention is given to learning models with general (equality and inequality) constraints or lacking differentiability. Such learning models are common in practice, especially in Intelligent Transportation Systems (ITS) that involve physical or domain knowledge as specific model constraints. Motivated by ITS applications, this paper formulates a model-target data poisoning attack as a bi-level optimization problem with a constrained lower-level problem, aiming to induce the model solution toward a target solution specified by the adversary by modifying the training data incrementally. As the gradient-based methods fail to solve this optimization problem, we propose to study the Lipschitz continuity property of the model solution, enabling us to calculate the semi-derivative, a one-sided directional derivative, of the solution over data. We leverage semi-derivative descent to solve the bi-level optimization problem, and establish the convergence conditions of the method to any attainable target model. The model and solution method are illustrated with a simulation of a poisoning attack on the lane change detection using SVM.},
keywords = {
Cybersecurity \sep Data poisoning attack \sep ITS \sep  Lipschitz continuity
}，
}

\begin{document}

\maketitle  

\Nomenclature

\begin{center}
\begin{tabular}{|P{.45\linewidth}|P{.45\linewidth}|}
\hline
$x$ & data, $x=[x_i],i=1,2,\cdots, n$ \\
\hline
$\Bar{x}$ & pristine data, $\Bar{x}=[\Bar{x}_i],i=1,2,\cdots,n $ \\
\hline
$x_i$ &  an individual data point, $i\in [1,2,\cdots ,n]$ \\
\hline
$x^k$ & poisoned data in the $k$th iteration \\
\hline
$x_i^k$ &  an individual poisoned data point in the $k$th iteration, $i\in [1,2,\cdots ,n]$ \\
\hline
$x^*$  & data that induces the target model \\
\hline
$\hat{x}$ & poisoned data that deviates the model (\ref{eqn:opt problem of ML}) the most\\
\hline
$\Delta x$  & feasible perturbation direction of $x$ \\
\hline
$\Delta x^k$  & the optimal perturbation direction of $x$ in the $k$th iteration \\
\hline
$\Delta x_i^k$  & the optimal perturbation direction of $x_i$ in the $k$th iteration\\
\hline
$S(x)$ & solution mapping from data to the model solution of model (\ref{eqn:opt problem of ML}), a single-valued mapping without specification\\
\hline
$y$   & optimization variable in learning model (\ref{eqn:opt problem of ML}) \\
\hline
$\hat{y}(x)$ & a solution of model(\ref{eqn:opt problem of ML}), i.e., $\hat{y}(x) = S(x)$ \\
\hline
$y^*$  &  target model, $y^*=S(x^*)$   \\
\hline
$\bar{y}$   & pristine solution to the learning model (\ref{eqn:opt problem of ML}), $\Bar{y}=S(\Bar{x})$ \\
\hline
$DG(x)(\Delta x)$ & semi-derivative of function $G(x)$ along the direction $\Delta x$\\
\hline 
$G'(x;\Delta x)$ & directional derivative of function $G$ at $x$ in the direction of $\Delta x$ \\
\hline

\end{tabular}
\end{center}

\section{Introduction}
\label{sec:intro}
Data poisoning attacks aim to add malicious noises, perturbations, or deviations to a training dataset to produce erroneous results if the dataset is used for training machine learning or deep learning (DL) models \citep{Akhtar2021Advances,Burkard2017Analysis,Mei2015Using, Wang2018Data}. As a well-recognized security risk, data poisoning attacks have been proposed and studied extensively, with attack methods developed for various learning models, such as linear regression (LR) \citep{Jagielski2018Manipulating}, support vector machines (SVM) \citep{Biggio2013Poisoning}, collaborative filtering \citep{Li2016Data}, neural networks \citep{Xie2022Universal}, etc. See Section \ref{sec:related} for detailed reviews on related works. With learning models continuing to serve as an important component in intelligent systems, data poisoning attacks are expected to increase in the future.  

The wide deployment of Internet of Things (IoT) devices in moving objects (humans, vehicles, etc.) has also promoted ``crowdsourcing'' as a popular means for data collection in many science and engineering fields. In Intelligent Transportation Systems (ITS), e.g., mobile sensing has become almost ubiquitous for data collection \citep{herrera2010evaluation,Ban2010Mobile}, which is a particular form of crowdsourcing so that users (human users, cars, bikes, scooters, etc.) contribute their individual travel/movement data (such as trajectories) to a central server or service (e.g.,  Google collects such data via mobile phone users who use Google Maps and opt-in to share data). Such data has revolutionized almost every aspect of ITS \citep{wang2024data}, which is, however, vulnerable to poisoning attacks due to its crowdsourcing nature: it is not hard for any user with malicious intent (or being hacked) to share poisoned data with the server to disrupt learning models that rely on the data. 

To better understand the vulnerability of the ITS system, data poisoning attacks have emerged these years (see section 2 for related works). The most prevalent attacks \citep{Mei2015Using,Biggio2013Poisoning,Papernot2016limitations} are to decrease the victim models' performance (e.g. prediction accuracy) on a given test dataset. The objective function is to maximize the deviation between the model solution and the pristine solution by injecting adversarial points into the training data. This adversarial goal, although natural and intuitive, can not allow us to customize the output/solution of the victim model after the attack, thereby failing to meet specific adversarial objectives. For example, during an attack on SVM, simply reducing the classification accuracy is insufficient to steer the decision hyperplane to a desired position, and thus cannot achieve specific classification objectives. This is the first motivation of this paper, i.e., to formally define a model-targeted data poisoning attack in ITS. For this, we need to find a target model that aligns with the attacker's objective, then induce the victim model towards the target by perturbing the train data.  
Secondly, ITS applications have unique features such that directly applying existing attack models may not be viable. One of such features is that ITS learning models often need to consider application-specific constraints. For example, a learning model to predict the vehicle queue at a signalized intersection often needs to consider the fact that the green time of the signal has a lower and upper bound \citep{Ban2011Real}, and two consecutive vehicles have to obey some physical car-following behavior. Such constraints can invalidate some widely-used attack models and algorithms, e.g., those based on gradients \citep{Mei2015Using}, because the model solution (e.g. model parameters) is no longer differentiable w.r.t. the data ({see section \ref{sec:gradient}}). This requires us to tailor and extend existing methods to better fit ITS applications. Therefore, our second motivation for this paper is to develop new methods/tools that rely on weaker conditions to define, analyze, and solve data poisoning attacks on ITS applications. Such methods and tools also generalize existing attack models and may be applicable to a broader set of learning models and their attack models in science/engineering applications.

This paper focuses on optimization-based learning models for ITS with possible equality and inequality constraints. The attack goal is to manipulate the models towards some adversarial, target models by perturbing the training data. We formulate the attack model as a bi-level optimization problem, with the lower level being the ITS learning model, while the upper level being the attack goal. To solve this bi-level problem, we need to measure the change of the model solution w.r.t. the data perturbation. This allows us to design data perturbation direction in each iteration to decrease the distance between the current solution and the target solution. Since the model solution is not differentiable w.r.t. the data due to the constraints of ITS models, we instead leverage the Lipschitz continuity to quantify the solution change w.r.t. the data perturbation. Under the assumption that the model solution is Lipschitz continuous w.r.t. the data, we can calculate the semi-derivative (also called directional derivative) of the solution over the data, which measures the solution's change in response to a given data perturbation. Note that the Lipschitz continuity assumption is weaker than the differentiability assumption used in most attack models. The semi-derivative-based method generalizes existing attack models and extends their application to a broader range of learning models. The resulting attack can guarantee the victim model achieves any target model with asymptotic convergence. We provide details of the Lipschitz-based analysis methods for poisoning attacks and theoretically analyze the convergence behaviors of the proposed attack models. We attack the lane change detection model using mobile trajectory data. We also compare the performances of the proposed model with state-of-the-art gradient-based methods. 
\section{Related Works}
\label{sec:related}
Although machine learning has achieved tremendous progress in traffic prediction \citep{Ban2011Real,wang2024transferability}, signal control \citep{guo2024network}, etc., it is not inherently designed with security considerations, making it vulnerable to adversarial attacks \citep{tian2022comprehensive}. Extensive studies have shown that data poisoning attacks can stealthily and significantly degrade model performance. The affected models range from classical machine learning algorithms such as support vector machines (SVMs) \citep{Biggio2013Poisoning} and regression models \citep{Jagielski2018Manipulating}, to more recent architectures including neural networks \citep{Papernot2016limitations} and large language models \citep{yao2024survey}. A representative example in the field is by  \citet{Biggio2013Poisoning}, who proposed a gradient ascent-based method to generate adversarial training examples for SVMs.  These attacks are typically implemented by altering input features or manipulating labels (e.g., through label flipping). Since label manipulation tends to be easier to detect, the feature modification approach has become more dominant in recent studies. For example, attackers can corrupt the SpamBayes learning models by injecting words likely to appear in legitimate emails into malicious ones, thus causing the model to misclassify spam as legitimate messages. 

Data poisoning attacks have been widely used to evaluate model robustness in domains such as malware classification \citep{suciu2018does}, natural language processing \citep{li2021hidden}, and web services \citep{xing2013take}. However, this threat has only recently begun to receive attention in the transportation domain.  Recently, researchers have developed attack models for traffic networks by poisoning only a small subset of spatiotemporal features \citep{wang2024transferability,liu2022practical,zhu2021adversarial}. These attacks take into account the vehicle speed range when constructing the perturbations. For more examples in the context of ITS, we refer the reader to \citep{wang2024datasurvey}.

Existing poisoning attacks may also be categorized by the objective functions: those that aim to maximally deviate the model solution from the pristine solution, and those that aim to deviate the poisoned model solution towards a target model. This paper focuses on the latter, i.e., model-targeted attacks.  While such attacks allow us to customize adversarial objectives, existing models are limited to unconstrained learning problems \citep{Suya2021Model-targeted}, as constraints invalidate gradient-based optimization methods for generating poisoned points. No matter which objective to consider, data poisoning attacks are often formulated as a bi-level optimization problem, where the upper-level problem is the adversarial objective, and the lower-level problem is the learning problem. Specifically, the attacker iteratively perturbs the data in a direction that helps approach the adversarial objective. To determine this direction, we need to capture the change in the adversarial objective corresponding to the data perturbations. This change is typically characterized using the gradient. However, when attacking learning problems with general constraints, such as ITS applications that involve physical and domain-knowledge-related constraints, it is well-known that the solution of the learning model is often continuous but not differentiable w.r.t. data changes \citep{luo1996mathematical}. Consequently, as the upper-level adversarial objective is defined using the solution of the learning problem, the objective function's change is not differentiable w.r.t. data changes (see section \ref{sec:gradient}). In such cases, one needs to turn to weaker conditions, like Lipschitz continuity \citep{Dontchev2009Implicit}, to capture the change of objective function and determine the data perturbation direction. This study introduces a semi-derivative descent algorithm to extend model-targeted attacks to convex models with constraints.

\section{Background and Problem Formulation}
\label{sec:background}
In this section, we first formulate the data poisoning attack as a bi-level optimization problem and then provide some comments on the gradient-based method to solve the attack model.

\subsection{Bi-level Optimization Model for Poisoning Attacks}

A learning model for ITS, which is also the victim model when data poisoning attack is concerned, can often be characterized as a constrained optimization problem.
\begin{equation}\label{eqn:opt problem of ML}
\begin{aligned}
& \underset{y}{\operatorname{min}} \quad g_0(x, y) \\
& s.t. \quad  g_i(x, y) \begin{cases}\leq 0 & \text { for } i \in[1, r], \\
=0 & \text { for } i \in[r+1, m].\end{cases}
\end{aligned} 
\end{equation}

Here, $y$ is the solution of the learning model (e.g., the coefficients of separating hyperplane in SVM or the weight vector of a neural network). $x \in R^d$ is the data used to train the model, which will be perturbed for poisoning attacks. We denote by $\mathcal{Y}(x)$ the feasible region of the model (\ref{eqn:opt problem of ML}). We define the value function by
\begin{equation}
    \mathcal{V}(x):=\inf \{g_0(x, y) \mid y \in \mathcal{Y}(x)\},
\end{equation}
and the solution mapping by
\begin{equation}
\mathcal{S}(x):=\{y \in \mathcal{Y}(x) \mid g_0(x, y)=\mathcal{V}(x)\}.
\end{equation}
This study focuses on data poisoning attacks given the attacker is capable of gaining full knowledge of the training data and model (i.e., white-box attack scenario). We aim to deviate the solution of the learning model \eqref{eqn:opt problem of ML}, i.e., $S(x)$, towards a pre-defined solution by minimally modifying the training data. We formulate the attack model (i.e., the process of finding a set of poisoning data points) as solving a bi-level optimization problem \citep{Mei2015Using}. Denote the adversarial goal by $G(x)$, the attack model can be written as below:

\begin{equation}\label{eq:2}
\begin{aligned}
& \underset{x}{\text{min}}
& & G(x)=|S(x)-y^*| \\
& s.t.
& & || x-\bar{x} \| \leq \delta, \\
&&& x \in X.
\end{aligned}
\end{equation}

Here $y^*$ is a predefined target model solution, i.e., there exists $x^* \in X$, such that $y^*\in S(x^*)$. For example, $y^*$ could be a classifier that mistakenly classifies trucks as bikes. The design of $y^*$ is based on the adversarial goal of attackers.  $G(x)$ is the attack objective, and $\bar{x}$ and $x$ are the pristine and poisoned data, respectively.  $\|x-\bar{x}\| \leq \delta$ states that the data perturbations shall be small to be stealthy, with $\delta$ application specific; $x \in X$ denotes other perturbation-related constraints required by the specific model. 

We can often easily design the target model solution $y^*$ to achieve our adversarial goal. However, we cannot obtain $x^*$, which is the training data that leads to $y^*$. By solving the model (\ref{eq:2}), we actually find a path from $\bar{x}$ to $x^*$, such that the model, when trained on $x^*$, can achieve $y^*$. 

The proposed attack model can be more flexible than existing data attacks \citep{Jagielski2018Manipulating}. By specifying the objective function as $G(x)=-\left|S(x)-S(\bar{x})\right|$, the objective is to induce the model to deviate from the original solution $S(\bar{x})$ as much as possible, see \citep{wang2024data}.

\subsection{Why Classical Gradient Descent Fails}
\label{sec:gradient}
The bi-level optimization problem (\ref{eq:2}) is to minimize the upper objective function $G(x)$  under the constraint that $S(x)$ is the solution of the lower optimization problem (\ref{eqn:opt problem of ML}). The gradient descent serves as a common, classical solution. Under very restrictive assumptions on the underlying learning model \eqref{eqn:opt problem of ML}, including smoothness and strong convexity of $g_0(x)$ and there is no constraint, the gradient of the upper objective function $G$ over data $x$ can be expressed using the chain rule.
\begin{equation}
\begin{aligned}
\nabla_x G(x)&= \nabla_y G(x,\hat{y}(x)) \nabla_x \hat{y}(x),
\end{aligned}
\end{equation}
where $\hat{y}(x) \triangleq \arg \min _{y \in \mathcal{Y}(x)} g_0(x, y) \in S(x)$ is the unique solution of model (\ref{eqn:opt problem of ML}) due to the strong convexity. When $\mathcal{Y}(x) = \mathcal{R}^{\operatorname{dim}(y)}$, by the (classical) implicit function theorem \citep{Dontchev2009Implicit} and optimality condition $\nabla_y g_0(x,\hat{y})=0$, $\nabla_x \hat{y}(x)$ can be expressed, which can help derive the gradient of the adversarial objective over data changes:
\begin{equation}
\begin{aligned}\label{eq:gradient of y}
\nabla_x \hat{y}(x)=-\nabla_{x \hat{y}}^2 g_0\left(x, \hat{y}(x)\right)\left[\nabla_{\hat{y} \hat{y}}^2 g_0\left(x, \hat{y}(x)\right)\right]^{-1} 
\end{aligned}
\end{equation}
\begin{equation}\label{eq:gradient}
\begin{aligned}
\nabla_x G(x)
=-\nabla_{x \hat{y}}^2 g_0\left(x, \hat{y}(x)\right)\left[\nabla_{\hat{y} \hat{y}}^2 g_0\left(x, \hat{y}(x)\right)\right]^{-1} \nabla_{\hat{y}} G\left(x, \hat{y}(x)\right)
\end{aligned}
\end{equation}

However, the gradient-based methods may not be applied if $\mathcal{Y}(x)$ is a region in Euclidean space, e.g., when the underlying learning model \eqref{eqn:opt problem of ML} has general equality/inequality constraints.  First, given the constraints in the lower problem (\ref{eqn:opt problem of ML}), the optimality condition $\nabla_y g_0(x,\hat{y})=0$ does not hold and should be replaced by the KKT condition, for which \eqref{eq:gradient of y} does not apply. Second, in this case, $\hat{y}(x)$ is often continuous, but not differentiable, with $x$ \citep{luo1996mathematical}, which leads to the nonexistence of gradient of $\hat{y}(x)$. To deal with attacks on general learning models, methods that can deal with the issues mentioned above are needed.

Figure \ref{fig:invalid gd} illustrates this via a simple bi-level optimization problem. The pristine data $x$ is at the origin and we try to perturb it along a feasible direction (towards right or left) to maximize the objective function $G(x,\hat{y}(x))$. Due to the constraint of the lower-level problem, $\hat{y}(x)$ is non-differentiable with $x$, which invalid formula (\ref{eq:gradient of y}) and (\ref{eq:gradient}).  The method we will propose next can properly identify the optimal perturbation direction to move the data (towards the right) with the help of the semi-derivative.

  \section{Lipschitz Continuity of Model Solution and Semi-derivative}
\begin{figure}[h!t]
\centering
\includegraphics[width=.6\linewidth]{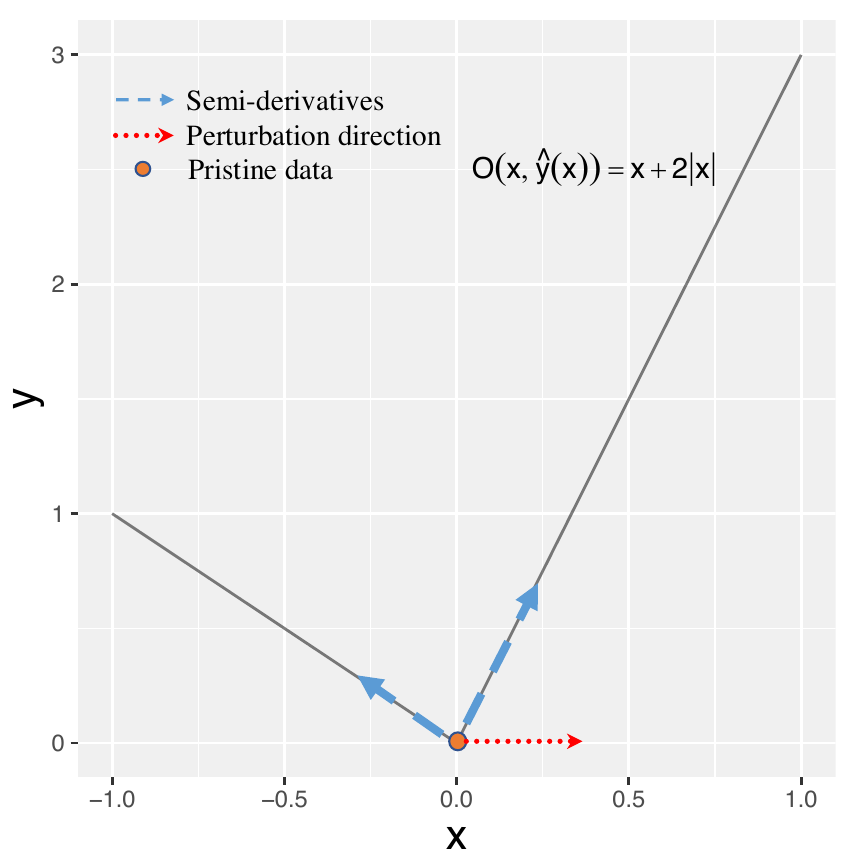}
\caption{Consider a bi-level optimization problem:$max_x G(x,\hat{y})=x+2\hat{y}$  $s.t.$  $ |x| \leq 1 $, $\hat{y} \in \operatorname{argmin}\{y: x+y\geq 0;y-x \geq 0\}$, the solution of lower problem is $\hat{y}=|x|$, and the upper objective function  $G(x,\hat{y}(x))= x+2|x|$, as shown above , is non-differentiable. Though the gradient is nonexistent, the semi-derivatives at pristine data w.r.t. all possible perturbation directions (right or left) do exist, which can be used to define and solve attacks.}
\label{fig:invalid gd}
\end{figure}

The analysis in Section \ref{sec:gradient} and the illustration in Figure \ref{fig:invalid gd} show that to properly analyze and solve the attack model, one needs to understand and characterize how the solution $\hat{y}(x)$ of the learning model \eqref{eqn:opt problem of ML} changes with data $x$, which is the classical sensitivity analysis of optimization problems \citep{Fiacco83}. The gradient is a natural choice, which, however, does not apply to the general learning model \eqref{eqn:opt problem of ML}. For this, we need to resort to definitions based on weaker conditions. Here we apply the concept of Lipschitz continuity and the semi-derivative. 

\begin{definition}[Lipschitz continuity]
\citep{Dontchev2009Implicit}) Assuming $x$ is in a $\delta$-neighborhood of $\bar{x}$, i.e., $\|x-\bar{x}\| \leq \delta$, and $K$ (the Lipschitz constant) is a finite value, we have the following definitions:

- If $\mathcal{S}(x)$ is \textbf{single-valued}: $\mathcal{S}(x)$ is Lipschitz continuous around $\bar{x}$ if $|\mathcal{S}(x)-\mathcal{S}(\bar{x})| \leq K|x-\bar{x}|$.

- If $S(x)$ is \textbf{set-valued}: $\mathcal{S}(x)$ is Lipschitz continuous around $\bar{x}$ if $S(x) \subset \mathcal{S}(\bar{x})+K|x-\bar{x}|\mathscr{B}$, with $\mathscr{B}$ a unit ball. 
\end{definition}

This paper focuses on the first case when (\ref{eqn:opt problem of ML}) has a locally unique solution. In this case, $S(x)$ contains only $\hat{y}(x)$. As the gradient of $\hat{y}(x)$ over $x$ may not exist, we define the semi-derivative of $\hat{y}(x)$ over $x$, as a counterpart of gradient, to capture how $\hat{y}(x)$ changes with data $x$.  

\begin{definition}[Semi-derivative \citep{Dontchev2009Implicit}]
     A function $f: \mathbb{R}^n \rightarrow \mathbb{R}^m$ is said to be semidifferentiable  at $\bar{x}$ if it has a first-order approximation at $\bar{x}$ of the form $h(x)=G(\bar{x})+\varphi(x-\bar{x})$ with $\varphi$ continuous and positively homogeneous.

The function $\varphi$ is called the semi-derivative of $f$ at $\bar{x}$ and denoted by $D G(\bar{x})$, so that $h(x)=G(\bar{x})+ D G(\bar{x})(x-\bar{x})$.
\end{definition}

``Lipschitz continuity'' of the solution $\hat{y}(x)$ (over data $x$) is the underlying condition we focus on in this study, which is weaker than the differentiability assumption needed for the gradient method. Under this condition, the semi-derivative of $\hat{y}(x)$ over $x$ exists and further semi-derivative is equivalent to the directional derivative 
\citep{Dontchev2009Implicit}. The Lipschitz constant can also be estimated by the upper bound of the semi-derivative,  which directly reflects how sensitive the solution $\hat{y}(x)$ changes with $x$. Figure \ref{fig:invalid gd} illustrates this semi-derivative when the solution mapping of the lower-level problem is nondifferentiable. 

Based on Lipschitz continuity and the semi-derivative, the proposed attack method is to select the data perturbation direction among all feasible perturbations each time to ensure $G(x)$, the distance between the current model solution and the target model solution, decreases as much as possible. Next, we show the derivation of semi-derivative for solution $\hat{y}(x)$.

\subsection{Semi-derivative of Solution Mapping } 
The semi-derivative of the solution can be derived using the general implicit function theorem for stationary points; see Theorem \ref{theorem:2G9} below, which can deal with the general learning model \eqref{eqn:opt problem of ML}. Denote $L(x,y,\lambda)$ the lagrangian function of (\ref{eqn:opt problem of ML}) with $\lambda$ the multipliers. $\bar{y},\bar{\lambda}$ are the (locally unique) solution of \eqref{eqn:opt problem of ML} at the pristine data $\bar{x}$. The derivation relies on the following \textbf{auxiliary problem} w.r.t. $\bar{x},\bar{y},\bar{\lambda}$:

\begin{equation}\label{eq:auxiliry}
  \begin{array}{cl}
\min \limits_{\Delta y} \bar{g}_0(\Delta y)-v^T(\Delta y), 
\text { s.t.} \quad \bar{g}_i(\Delta y)+\mu_i \begin{cases} = 0 & \text { for } i \in I \backslash I_0 ; \\
\leq 0 & \text { for }i \in I_0.
\end{cases}
\end{array}  
\end{equation}
With $I, I_0$ two index sets defined as,
\begin{equation*}
\begin{aligned}
I =\left\{i \in[1, m] \mid g_i(\bar{x}, \bar{y})=0\right\}, 
I_0 =\left\{i \in[1, r] \mid g_i(\bar{x}, \bar{y})=0 \text { and } \bar{\lambda}_i=0\right\}.
\end{aligned}
\end{equation*}
And $\bar{g}_0$ is the second-order approximation of the objective function and $\bar{g}_i$ is the first-order approximation of the constraint, making the auxiliary problem a quadratic program with linear constraints.
\begin{equation*}
\begin{aligned}
\bar{g}_0(\Delta y) & =L(\bar{x}, \bar{y}, \bar{\lambda})+\left\langle\nabla_y L(\bar{x}, \bar{y}, \bar{\lambda}), \Delta y\right\rangle +\frac{1}{2}\left\langle\Delta y, \nabla_{y y}^2 L(\bar{x}, \bar{y}, \bar{\lambda}) \Delta y\right\rangle. \\
\bar{g}_i(\Delta y) &=g_i(\bar{x}, \bar{y})+ \left\langle\nabla_y g_i(\bar{x}, \bar{y}, \bar{\lambda}), \Delta y\right\rangle.
\end{aligned}
\end{equation*}

Note that $\Delta y$ denotes the perturbation to the solution $\bar{y}$; $v$ and $\mu$ are two parameters of the auxiliary problem that need to be specified when defining the problem.

The auxiliary problem  (\ref{eq:auxiliry}) serves as a bridge that helps connect the data perturbation in $x$ with changes in the solution of \eqref{eqn:opt problem of ML}. In particular, under certain conditions, the solution of the auxiliary problem helps characterize the semi-derivative of $\hat{y}$ over $x$. To formally derive this, note that when \eqref{eqn:opt problem of ML} has a single solution, it is $\hat{y}(\bar{x})$ for the input data $\bar{x}$. Then the semi-derivative of $\hat{y}(\bar{x})$ along the direction $\Delta x$, denoted as $D \hat{y}(\bar{x})(\Delta x)$, can be calculated via the general implicit function theorem, as stated in Theorem \ref{theorem:2G9} below. It directly follows Theorem $2 \mathrm{G} .9$ in \citep{Dontchev2009Implicit} and thus the proof is omitted.

\begin{theorem}
Suppose the optimization model expressed in (\ref{eqn:opt problem of ML}) with twice continuously differentiable function $g_i$ and the following conditions hold:
 \begin{itemize}
     \item [i)]The gradients $\nabla_y g_i(\bar{x}, \bar{y})$ for $i \in I$ are linearly independent.
     \item[ii)] $\left\langle\Delta y, \nabla_{y y}^2 L(\bar{x}, \bar{y}, \bar{\lambda}) \Delta y\right\rangle>0$ for every nonzero $\Delta y \in M^{+}$, a subspace defined as:
$$
M^{+}=\left\{\Delta y \in R^t \mid \Delta y \perp \nabla_y g_i(\bar{x}, \bar{y}) \text { for all } i \in I \backslash I_0\right\}.
$$
 \end{itemize}
Then, the solution set of \eqref{eqn:opt problem of ML}, $S(x)$, has a Lipschitz continuous single-valued localization $\hat{y}$ around $\bar{x}$, and this localization $\hat{y}$ is semi-differentiable at $\bar{x}$ with the semi-derivative given by
\begin{equation}\label{eq:theorem1}
\begin{aligned}
& D \hat{y}(\bar{x})(\Delta x)=\bar{s}(-B \Delta x), \\
\text{with} \quad & B=\left(\begin{array}{c}
\nabla_{y x}^2 L(\bar{x}, \bar{y}, \bar{\lambda}) \\
-\nabla_x g_i(\bar{x}, \bar{y}) \text { for } i \in[1, m]
\end{array}\right).
\end{aligned}
\end{equation}

\label{theorem:2G9}
\end{theorem} 

Here $\bar{s}(z)$ denotes the solution of the auxiliary problem, $z=-B\Delta x$ is the input that defines the problem. $B$ represents the Hessian matrix of the Lagrangian and the gradient of the constraints of the learning problem (\ref{eqn:opt problem of ML}) evaluated at $\bar{x}$. Their multiplications with $\Delta x$, as shown in (\ref{eq:theorem1}), produce the perturbations to the auxiliary objective ($v$) and constraints ($\mu$), which are required as the input to define the auxiliary problem. In other words, the calculation of the semi-derivative of the solution of (\ref{eqn:opt problem of ML}) at $\bar{x}$ is as follows: the first step is to solve (\ref{eqn:opt problem of ML}) and obtain (unique) solution $\bar{y}, \bar{\lambda}$, and matrix $B$; the second step is to construct the auxiliary problem at $\bar{x}$; the third step is to solve the auxiliary problem and obtain its solution to calculate the semi-derivative as shown in \eqref{eq:theorem1}. 

\subsection{Semi-derivative Descent}
This section presents a semi-derivative descent algorithm for solving the bi-level attack model (\ref{eq:2}). This approach is analogous to gradient descent, but replaces the gradient with semiderivative to overcome the non-differentiability. A pseudocode is provided in Appendix B.

Using Theorem \ref{theorem:2G9}, we can compute the semi-derivative $DG(x)( \Delta x)$, the change of $G(x)$ when we perturb the data $x$ along direction $\Delta x$. By the chain rule, $DG(x)(\Delta x)$ can be expressed as:
\begin{equation}\label{eq:8}
\begin{aligned}
DG(x)(\Delta x) & =\nabla_{\hat{y}} G(x) \cdot D\hat{y}(x)(\Delta x),
\end{aligned}
\end{equation}
where the derivation of $\nabla_{\hat{y}} G$ takes advantage of the fact that $G$ is often an explicit function of $\hat{y}(x)$. $D\hat{y}(x)(\Delta x)$ can be derived by (\ref{eq:theorem1}), which could often be done by solving a set of equations. Note also that the semi-derivative (\ref{eq:8}) will be equivalent to the gradient (\ref{eq:gradient}) if $\hat{y}(x)$ is differential with $x$ and there is no constraint in the lower problem (\ref{eqn:opt problem of ML}). In the $k$th iteration, when the data has been modified as $x^k$ during the previous steps, we find the $\Delta x^k$ s.t. 
\begin{equation}
    \Delta x^k=\underset{\Delta x}{\operatorname{argmin}} D G\left(x^k\right)(\Delta x).
\end{equation}
We then perturb the data $x^k$ along the direction of $\Delta x^k$ with step size $\eta^k$. Note that we always adopt $\Delta x^k$ with the most negative semi-derivative to ensure $G(x)$ will decrease after $x^k$ moves along $\Delta x^k$. If $\underset{\Delta x}{\operatorname{argmin}} \, D G\left(x^k\right)(\Delta x) \geq 0$, which indicates that $G(x)$ cannot decrease along any feasible perturbation direction of $x^k$, $x^k$ achieves the local/global minima. 

\subsection{Convergence of Semi-derivative Descent}

Before illustrating the convergence result, this section first formally introduces the definitions, assumptions, and optimality conditions that serve the convergence of semi-derivative descent. 
\begin{definition}[Attainable target model]\label{def:Attainable}
    We say that $y^*$ is an attainable target model if $y^*$ is a feasible solution of model (1), i.e.,
\begin{equation}
\exists x^* \in \Omega \text {, s.t. } y^* \in S\left(x^*\right) \text {. }
\end{equation}
\end{definition}
Here, $S$ stands for solution mapping from $X$ to $Y$, and $\Omega=X \cap B(x, \delta) \subset \mathbb{R}^n$ denotes the feasible region with respect to $x$, with $B(x, \delta)$ denoting possible perturbation region of $x$. Following Definition \ref{def:Attainable}, $G\left(x^*\right)=|S(x^*)-y^*|=0$, i.e., the poisoned data $x^*$ that induces the (poisoned) solution mapping $\hat{y}$ converge to $y^*$ is in fact the minimizer of the $G(x)$ over $x \in \Omega$:
\begin{equation}
x^*=\underset{x \in \Omega}{\operatorname{argmin}} G(x).
\end{equation}

\begin{definition}[Feasible attack direction]\label{def:Feasible attack direction}
    $\Delta x$ is a feasible attack direction if $x+\epsilon \Delta x \in \Omega$ for some $\epsilon>0$.
The set of feasible attack directions is denoted by $Q_x$ :
$$
Q_x=\left\{\Delta x \in R^{\operatorname{dim}(x)} \mid x+\epsilon \Delta x \in \Omega\right\}.
$$
\end{definition} 

\textbf{Remark:} \textit{In the \textit{k}th iteration, \textit{Algorithm 1} perturbs the data $x^k$ in the direction $\Delta x^k$ such that the semi-derivative of $G(x)$ is most negative at $x^k$ along $\Delta x^k$ among all feasible attack directions, i.e., 
\begin{equation}
    DG(x^k)(\Delta x^k)= \mathop{argmin} \limits_{\Delta x \in Q_x \atop ||\Delta x||=1} DG(x^k)(\Delta x). \label{eq:p}
\end{equation}}
\begin{corollary}[Optimality condition] \label{corollary:optimality}
Let $G(x)$ be a semi-differentiable function with respect to $x$, then
\begin{equation}
\forall \Delta x \in Q_x,\|\Delta x\|=1, D G\left(x^*\right)(\Delta x) \geq 0 \quad \Rightarrow x^* \in \underset{\Omega}{\operatorname{argmin}} G.
\end{equation}

\end{corollary}
\textit{proof: see Appendix A.2.}

\textbf{Remark}. Note that when $x=x^* $ s.t. $ DG(x^*)(\Delta x) \geq 0 $ for any feasible attack direction $\Delta x$, $G(x^*)$ reaches the minimum over $x \in \Omega$. For the optimal perturbation direction in the $kth$ iteration, $\Delta x ^k=argmin_{\Delta x} DG(x^k)(\Delta x)$, Corollary \ref{corollary:optimality} shows that $DG(x^k)(\Delta x^k)$ is always negative before the algorithm converges to the target model, i.e., before $x$ reaches the global minimizer $x^*$ of $G(x)$.

Theorem \ref{theorem:1} below establishes the convergence of semi-derivative descent to any attainable target model $y^*$. We begin by stating the necessary assumptions, followed by the presentation of Theorem \ref{theorem:1}.
 
\textbf{Assumption 1.} The function $G(x)$ in model (\ref{eq:2}) is $\sigma$-strongly convex (see definition \ref{Strong convexity})  and twice directional differentiable (see definition \ref{Twice directional differentiable}). The absolute value of the second-order directional derivative of $G$ is bounded by $L$. We can assume $L>\sigma$ because $L$ is an upper bound.

\textbf{Assumption 2.} Mapping $S$ is uniformly Lipschitz continuous.
\textbf{Assumption 3.} $\Omega \subset R^n$ is compact.
\begin{theorem}
Let $\eta^k=-DG(x^k)(\Delta x^k)/L$. Assume that $G$ is strongly convex with modulus $\sigma$, and the norm of the second-order directional derivative of  $G(x)$ is bounded by $L$, then for $k>0$:
\begin{equation}
    G(x^*)-G(x^k)\leq (1-\frac{\sigma}{L})^{k}(G(x^*)-G(\bar{x}))). 
\label{9}
\end{equation}
\label{theorem:1}
\end{theorem}
\textit{proof: see Appendix A.2.}

Here $\Bar{x}$ is the pristine data, $x^k$ is the poisoned data after the \textit{ kth} iteration, $x^*$ is the poisoned data that produces the attainable target model, i.e., $y^*\in S(x^*)$. 
\section{Implementing Semi-derivative Descent for Data Poisoning Attack} 
We apply the semi-derivative descent to solve the model (\ref{eq:2}), to implement the data poisoning attack. We perturb only one data point each time to ensure that the data modification is imperceptible. The basic idea of the proposed attack is given by Algorithm \ref{alg:1}, see Appendix B. The attack algorithm sequentially applies perturbations to the dataset. In each iteration, we select one data point with the maximum effect (i.e., semi-derivative with the largest absolute value) and perturb it along the most negative semi-derivative; the learning model defined in (\ref{eqn:opt problem of ML}) is then resolved on the perturbed data after each iteration until convergence conditions are met. 
\subsection{Experiments}

\textbf{(Experimental Setup)} We implement the proposed attack methods on an SVM-based lane change detection using part of the vehicle trajectory data from the NGSIM US-101 dataset. The details of this model can be found in Appendix C. The SVM-based lane change classifier utilizes lateral velocity (m/s) and space speed (m) as input features to predict whether a vehicle will change lanes. The lateral velocity is calculated by differentiating the vehicle's lateral position over time. The binary lane change label is determined by observing transitions in the Lane ID field over time. 

The weight vector of the separating hyperplane is denoted by $w=[w_1, w_2]=[-7.64,-13.34]$, where $w_1, w_2$ represent the coefficient of normalized lateral velocity and headway space. Our adversarial goal is to induce the model solution ${w}=({w_1},{w_2})$ towards a target solution $\boldsymbol{w}^*=\left(w_1^*, w_2^*\right)$, where $\left|w_1^*-w_2^*\right|=0$, aiming to equalize the weights of different features as much as possible. The objective function of (\ref{eq:2}) is specific to minimize $\left|w_1-w_2\right|$. Using  $x_i^{v}$,$x_i^{h}$, and $y_i$ to denote the lateral speed, headway space, and label of vehicle $i$, the attack model is given by 
\begin{equation}\label{eq:attacksvm}
\begin{aligned}
\min_{\boldsymbol{x}} \quad & \|\boldsymbol{w}_1 - \boldsymbol{w}_2\|_2^2 \\
\text{s.t.} \quad & |\boldsymbol{x} - \bar{\boldsymbol{x}}| \leq \delta, \\
&\rho_1 \leq x_i^{v} \leq \rho_2, \quad \rho_3 \leq x_i^{h} \leq \rho_4, \\
& (\boldsymbol{w}_1, \boldsymbol{w}_2) \in \arg\min_{\boldsymbol{w}, b} \ \frac{1}{2} \|\boldsymbol{w}\|^2 + C \sum_{i=1}^n \xi_i, \\
& \quad \text{s.t. } \quad y_i(\boldsymbol{w}^\top \bar{x}_i + b) \geq 1 - \xi_i,\quad \xi_i \geq 0, \quad \forall i
\end{aligned}
\end{equation}
where $\bar{x}_i=(\bar{x}_i^v,\bar{x}_i^h)$ is the original (unperturbed) feature vector. $\rho_1$-$\rho_4$ define the domain-specific bounds on lateral speed and headway to ensure physical feasibility.

\begin{figure}[h]
  \centering
  \includegraphics[width=0.8\textwidth]{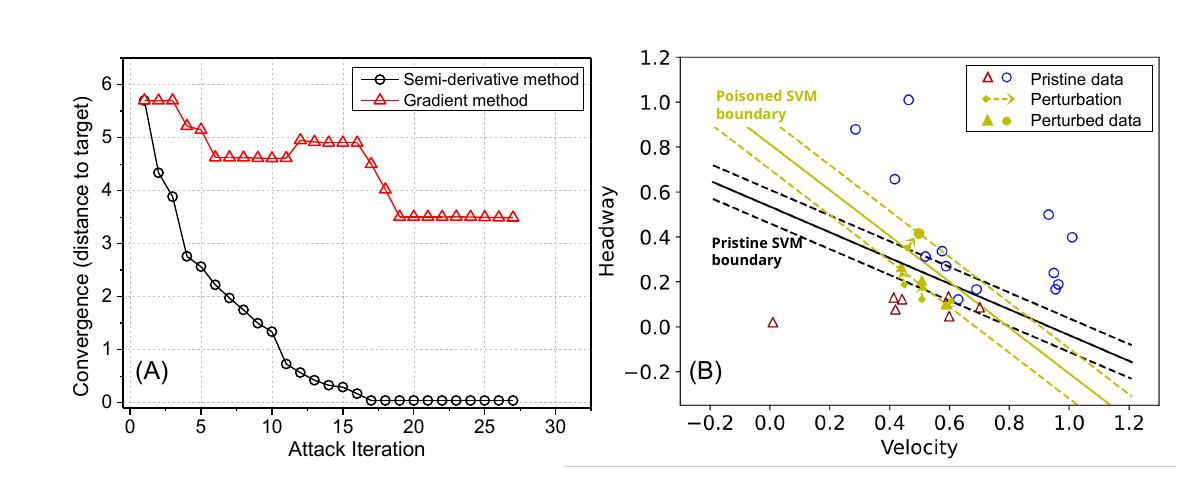}
  \caption{ A) Comparing convergence of the proposed attack with the gradient method. B) Poisoned model and data following the proposed attack.}
  \label{fig:fig2}
\end{figure}

We derive the semi-derivative of the bi-level optimization problem (\ref{eq:attacksvm}) to implement the semi-derivative descent.  We also compare the results with the gradient-based attack methods \citep{Biggio2013Poisoning,Jagielski2018Manipulating}. Figure \ref{fig:fig2} compares the convergence to the adversarial objective following the semi-derivative and the gradient method. We can observe that the semi-derivative attack steadily rotates the hyperplane and consequently converges to the target model. After the convergence, four data points are slightly perturbed and the SVM’s decision boundary is significantly shifted (Figure \ref{fig:fig2}B), making the two features equally important. In contrast, the gradient attack does not seem to converge toward the target model reliably (Figure \ref{fig:fig2}A). The gradient attack not only converges at a much slower rate than the semi-derivative method but also stops early after only a few iterations. Readers may see more experimental results of the proposed method, applied to both SVM and regression models in \cite{wang2024data}.
\section {Concluding Remarks}
This work formulates data poisoning attacks against machine learning (ML) models as sensitivity analysis of optimization problems over data perturbations (attacks) and studies the Lipschitz continuity property of the solution of attacked ML models. Based on the Lipschitz continuity, we propose methods to calculate the semi-derivative of the solution over data perturbations, which enables us to develop attack models to fit a broader spectrum of learning tasks, especially those with general constraints. The proposed method has proved convergence, which is applied to and tested using an ML-based ITS model. The results show that the proposed method can outperform the state-of-the-art gradient method in terms of flexibility and the effectiveness of attacks. We demonstrate that without considering the role of constraints, the gradient methods would be slow in moving toward the adversarial objectives or even fail to advance when the constraints become active. The study has multiple implications, including understanding the vulnerability of ML-based ITS models and developing effective defense methods.

This study assumes that the model (\ref{eqn:opt problem of ML}) has a unique solution $\hat{y}(x)$, which is common under the ML scenario. When the model (\ref{eqn:opt problem of ML}) represents a deep learning ITS model, such as a neural network, the solution mapping $S(x)$ may not be a singleton. Our future research aims to extend the Lipschitz analysis method and the calculation of semi-derivatives to set-valued maps, i.e., when the underlying learning model has a set of solutions instead of a single solution. In the set-valued case, the change of $S(x)$ in response to the data perturbation along a given direction $\mu$  can be measured following $\underset{\substack{\tau \rightarrow 0}}{\limsup } \frac{S\left(\bar{x}+\tau \mu\right)-S(\bar{x})}{}$. This limit measures the expansion and contraction of the model solution set $S(x)$ w.r.t. the data perturbation, which generalizes the semi-derivative we used in this paper, allowing us to design the data poisoning attack for DL models. We are working on this in our current research and results may be reported in subsequent papers.

\section*{Acknowledgement}
The research of the first, fourth, and fifth authors is partially supported by the National Science Foundation (NSF) grant CMMI-2326340. The research of the third author is partially supported by the National Science Foundation (NSF) grant CMMI-2326341. The work of the second author was conducted while a Ph.D. student at the University of Washington. Any opinions, findings, conclusions, recommendations expressed in this paper are those of the authors and do not necessarily reflect the views of NSF.
\newpage
\appendix
\appendix
\renewcommand{\thesection}{Appendix~\Alph{section}}
\section*{Appendix A: Convergence for Target Poisoning Attack}\label{sec:appendixa}
\addcontentsline{toc}{section}{Appendix A: Convergence for Target Poisoning Attack}
\label{appendixA}

This appendix shows that, under certain conditions, semi-derivative descent can converge to an attainable target model $\mathbf{y}^*$ when it is applied to the model (\ref{eq:2}). 
We first give a few definitions and the main assumptions to establish the convergence results. In the appendix, we have renumbered and restated all necessary definitions, theorems, and corollaries.

\subsection*{Appendix A.1. Definition and Assumptions}
\addcontentsline{toc}{section}{Appendix A.1.Definition and Assumptions}
\label{appendixA.1}

In this paper, we assume that $\hat{y}(x) \triangleq \arg \min _{y \in \mathcal{Y}(x)} g_0(x, y) \in S(x) $  is the unique solution of model (\ref{eqn:opt problem of ML}). Given the objective function of the model (\ref{eq:2}), $G(x)=|S(x)-{y}^*|$, where ${y}^*=S({x^*})$ is the attainable target model, to convenient our discussion, model (\ref{eq:2}) can be expressed as a minimization problem: 
\begin{equation}\label{eq:new-bilevl}
\begin{aligned}
\min_{x} \quad & G(x)= |\hat{y}(x)-y^*|\\
\textrm{s.t.} \quad & x \in \Omega, \\
  & \hat{y} \in S(x),   \\
\end{aligned}
\end{equation}
where $|\hat{y}(x)-y^*|$ is the squared distance towards the target model $y^*$, and without specification, $|\cdot|$ denots the squared $L_2$ norm. Following the definition \ref{def:Attainable}, the poisoned data $x^*$ that induces the (poisoned) solution mapping $\hat{y}$ converge to $y^*$ is in fact the  minimizer of the model (\ref{eq:new-bilevl}) over $x\in\Omega$:
\begin{equation}
    x^*=\mathop{argmin} \limits_{x \in \Omega} G(x).
\end{equation}

\begin{definition}[Directional differentiability]
Let $X \subset \mathbb{R}^n$ be an open set,  $x^0 \in X$ and  $G: X \longrightarrow \mathbb{R}$. If limit
$$
\lim _{t \rightarrow 0^{+}} \frac{G\left(x^0+t v\right)-G\left(x^0\right)}{t}
$$
exists (finite or not), the limit is the (one-sided) directional derivative of $G$ at $x^o$ in the direction of $v$. If the directional derivative exists for every $v$, $f$ is said to be directionally differentiable at $x^0$.
\end{definition}

\textbf{Remark.} \textit{Under the condition of Lipschitz continuous, semi-differentiability is equivalent to directional differentiability :
\begin{equation}
\lim _{t \rightarrow 0^{+}} \frac{G(\bar{x}+t \Delta x)-G(\bar{x})}{t} := G^{\prime}(\Bar{x};\Delta x)=  DG(\bar{x})(\Delta x). \label{eq:d}
\end{equation}
Note that $DG(\bar{x},\Delta x)$ is positively homogeneous:
\begin{equation}
    DG(\bar{x})(\alpha\Delta x)=\alpha DG(\bar{x})(\Delta x), \quad \text{whenever} \quad \alpha >0.
\end{equation}}

\begin{definition}[Twice directional differentiable]\label{Twice directional differentiable}  We say a function $G: X \rightarrow \mathbb{R}$ is twice directional differentiable at $u \in X$ if the limit

\begin{equation}
G^{\prime \prime}(u ; \varphi, \psi)=\lim _{\theta \rightarrow 0}\left(G^{\prime}(u+\theta \psi; \varphi)-G^{\prime}(u; \varphi)\right) / \theta
\end{equation}
exists for each $\varphi, \psi \in X$. A function $G: X \rightarrow \mathbb{R}$ has a second-order directional Taylor expansion at $x^0$ if
\begin{equation}\label{eq:equivalent}
 G\left(x^0+t \phi\right)=G\left(x^0\right)+tG^{\prime}(x^0;\phi)+\frac{t^2}{2}G^{\prime\prime}(x^0+t^0\phi;\phi,\phi),
\end{equation}
where $t^0\in [0,t]$.
\end{definition}

\begin{definition}[Strong convexity]\label{Strong convexity} $f$ is $\sigma$-strongly convex with respect to the Euclidean norm if there exists a modulus of convexity $\sigma>0$ such that
$$
G(y) \geq G(x)+G(x;y-x)+\frac{\sigma}{2}\|y-x\|_2^2, \quad \text { for all } x, y .
$$
(Henceforth, we use $\|\cdot\|$ to denote the Euclidean norm $\|\cdot\|_2$, unless otherwise specified.)
\end{definition}
\begin{assumption}
The function $G(x)$ in (\ref{eq:new-bilevl}) is $\sigma-$strongly convex and twice directional differentiable. The absolute value of the second-order directional derivative of $G$ is bounded by $L$. We can assume $L>\sigma$ because L is an upper bound. 
\end{assumption}
\begin{assumption}
\label{eq:assump2}
Mapping $S$ is uniformly Lipschitz continuous. 
    
\end{assumption}
\begin{assumption}
    
\label{eq:assump3} $\Omega \subset R^n$ is compact.  

\end{assumption}
\textbf{Remark.} \textit{Assumption 3 above is an assumption often applied for similar analysis, which should hold for most of the attack models. Assumption 2 seeks the Lipschitz continuity of the solution mapping $S$ of the model (\ref{eqn:opt problem of ML}), which we focus on in this study. Assumption 1 concerns the objective function of the attack model, which is satisfied if $L_2$ norm is used in (\ref{eq:new-bilevl}).}
 
\subsection*{Appendix A.2. Convergence Analysis}
\addcontentsline{toc}{section}{Appendix A.2. Convergence Analysis}
\label{proof of theorem and collary}

\textbf{\textit{Proof of corollary 1 :}} As an application of (\ref{eq:d}), $DG(x^*)(\Delta x) \geq 0$ implies that
\begin{equation}
   \lim _{t \rightarrow 0} (G(x^*+t\Delta x)-G(x^*))/t \geq 0.
\end{equation}
Hence, there exists $\delta >0$, such that for $\forall t \in (0,\delta)$
\begin{equation}
\label{eq:temp1}
G(x^*+t\Delta x)-G(x^*)\geq 0,
\end{equation}
where $ \Delta x \in Q_x$, and $||\Delta x||=1$.
For any $x \in \mathscr{B}(x^*,\delta)\bigcap \Omega$, let $q=\frac{x-x^*}{||x-x^*||}$. Since $x \in \mathscr{B}(x^*,\delta)$, we have$||x-x^*||\leq \delta$;  using (\ref{eq:temp1}), we obtain
\begin{equation}
    G(x)=G(x^*+||x-x^*||\frac{x-x^*}{||x-x^*||}) \geq G(x^*).
\end{equation}
Following the condition that $f$ is a convex function, $x^*$ is a minimizer of $f$ over $x \in \Omega$. $\Box$

\textbf{\textit{Proof of theorem 2}}. By applying Taylor expansion in the sense of directional derivative and using the assumption that the second-order directional derivative is bounded by $L$, we have
\begin{equation}
\label{eq:bound}
\begin{aligned}
    G(x^{k+1})&=G(x^k+\eta^k\Delta x^k),\\
    &= G(x^k)+\eta^kG^{\prime}(x^k;\Delta x^k)\\&+\frac{(\eta^k)^{2}}{2}f^{\prime\prime}(x^k+\eta_0\Delta x^k;\Delta x^k,\Delta x^k),\\
    &\leq G(x^k)-\frac{1}{2L}[G^{\prime}(x^k;\Delta x^k)]^2.
\end{aligned}
\end{equation}
Note that $\eta^k=-DG(x^k)(\Delta x^k)/L= G^{\prime}(x^k;\Delta x^k)/L$, and $\eta_0 \in [0,\eta^k]$. \\

By $\sigma$-strong convexity of $f$, for  any $z \in \Omega$, we obtain
\begin{align}
   G(z) \geq G(x^k)+G^{\prime}(x^k;z-x^k)+\frac{\sigma}{2}||z-x^k||^2_2.
\end{align}
Following the positively homogeneous and (\ref{eq:p}), we then have
\begin{align}
    G(z) &\geq G(x^k)+||z-x^k||G^{\prime}(x^k;\frac{z-x^k}{||z-x^k||})+\frac{\sigma}{2}||z-x^k||^2_2,\\   
    &\geq   G(x^k)+||z-x^k||G^{\prime}(x^k;\Delta x^k)+\frac{\sigma}{2}||z-x^k||^2_2.
\end{align}
Minimizing the left side with respect to $z$ over $z\in \Omega$ and minimizing the right side with respect to $z$ over $z\in R^n$, we have

\begin{align}
    G(x^*)&=\mathop{min} \limits_{z \in \Omega} G(z), \\
               &\geq \mathop{min} \limits_{z \in \Omega}  G(x^k)+||z-x^k||G^{\prime}(x^k;\Delta x^k)+\frac{\sigma}{2}||y-x^k||^2_2,\\
                &\geq \mathop{min} \limits_{z \in R^n} G(x^k)+||z-x^k||G^{\prime}(x^k;\Delta x^k)+\frac{\sigma}{2}||z-x^k||^2_2. \label{eq:12}
\end{align}

Let $M=||z-x^k||$, the right side of (\ref{eq:12}) achieves the minimal value when $M=-G^{\prime}(x^k;\Delta x^k)/\sigma$.  \\
Note that $G^{\prime}(x^k;\Delta x^k)$ is always negative until  $x^k$ converges to the minimizer of $G(x)$. In (\ref{eq:12}), substituting $||y-x^k||$ with $-G^{\prime}(x^k;\Delta x^k)/\sigma$ leads to 
\begin{align}
      G(x^*)&=\mathop{min} \limits_{y \in \Omega} G(y) \geq G(x^k)-\frac{1}{2\sigma}[G^{\prime}(x^k;\Delta x^k)]^2. \label{eq:b}
\end{align}
If $x^k$ is the minimizer of $G(x)$ over $\Omega$, then $G(x^k)=G(x^*)$, (\ref{eq:b}) still holds. 

Following (\ref{eq:bound}) and (\ref{eq:b}), we have 
\begin{align}
    G(x^{k+1})-G(x^*) \leq (1-\frac{\sigma}{L})(G(x^k)-G(x^*)).\label{eq:c}
\end{align}
Then, the convergence by (\ref{9}) in \textbf{Theorem \ref{theorem:1}} can be obtained via recursively applying (\ref{eq:c}). $\Box$
\section*{Appendix B: Model-targeted Data Poisoning Attack} 
\addcontentsline{toc}{section}{Appendix B: Model-targeted Data Poisoning Attack} 
\label{appendixB}



\begin{algorithm}[H]\label{alg:1}
\renewcommand{\algorithmicrequire}{\textbf{Input:}}
	\renewcommand{\algorithmicensure}{\textbf{Output:}}
    \caption{Semi-derivative-based Data Poisoning Attack Algorithm}
    \label{alg:1}
    \begin{algorithmic}[1]
    \REQUIRE Clean data $\bar{x}=\left[\bar{x}_i\right]\ (i\in\left[1,\ldots,n\right])$ and a victim model. \eqref{eqn:opt problem of ML}\\             \REPEAT
        \STATE For each data point $x_i^k$ at iteration $k$, perturb the individual data point along a random unit perturbation $\delta_i$, and evaluate the change of $G(x^k)$ in response to this perturbation by its semi-derivative $DG\left(x^k\right)\left(\delta_i \right).$

        \STATE Select the most sensitive data points, with the index denoted by $p \in \left[1,\ldots,n\right]$. 
\begin{equation*}
\begin{gathered}
p=\mathop{argmax}\limits_{i} \|DG\left(x^k\right)\left(\delta_i\right)\|, i=1,2,\cdots,n. 
\end{gathered}
\end{equation*}
\vspace{-0.1in}
        \STATE For data point $x_p^k$, derive its perturbation direction that leads $G(x)$ decrease most.  
\begin{equation*}
\Delta x_p^k=\underset{\Delta x_p}{\operatorname{argmin}} D G\left(x^k\right)(\Delta x_p).
\end{equation*}

\vspace{-0.1in}
        \STATE Perturb the target point $x_p$ following direction $\Delta x_p^k$ by a specific step size $\eta$.
\begin{equation*}
x_p^{k+1}=x_p^k+\eta^k \Delta x_p^k.
\end{equation*}
\vspace{-0.1in}
        \STATE Re-learn the victim model on the poisoned data and evaluate the attack performance.
        \UNTIL{Improvement in the attack performance is sufficiently small or constraints are violated.}
    \ENSURE Poisoned data and corresponding solution of victim model.
    \end{algorithmic}
\end{algorithm}

\section*{Appendix C: SVM based lane change detection} 
\addcontentsline{toc}{section}{Appendix C: SVM based lane change detection}
\label{appendixC}
This section provides the details of the SVM-based lane change detection model, which serves as the victim model in our experiments.

\textbf{Data Description.} The dataset used in this study is extracted from the NGSIM US-101 dataset, which provides high-resolution vehicle trajectory data captured at 0.1-second intervals. We compute the lateral velocity of each vehicle by differentiating its lateral position with respect to time (i.e., dividing the change in lateral position by 0.1 seconds). After preprocessing, each data point consists of the following characteristics: lateral velocity, the distance from the preceding vehicle in the same lane (space headway), and a binary label indicating lane change behavior. A value of 1 represents a lane change, while -1 indicates no lane change. We denote the normalized lateral velocity and normalized headway as \( x_i^v \) and \( x_i^h \), respectively.

\textbf{Model and Results.} We implement a binary classifier using a linear Support Vector Machine (SVM) to predict whether a vehicle will perform a lane change. The model is trained using the \texttt{SVC} class from the \texttt{scikit-learn} package with a linear kernel. Feature standardization is performed using \texttt{StandardScaler}. After training, the separating hyperplane, support vectors, and margin boundaries are visualized in the figure \ref{fig:fig2}. We report the parameters of the model below.

\begin{table}[htbp]
\centering
\caption{Classification performance and learned SVM parameters}
\begin{tabular}{lccc}
\toprule
Metric & Precision & Recall & F1-Score \\
\midrule
SVM Result & 1.00 & 0.93 & 0.96 \\
\midrule
Weight vector \( \boldsymbol{w} \) & \multicolumn{3}{c}{\( [-7.64,\ -13.34] \)} \\
\bottomrule
\end{tabular}
\label{tab:svm_result_clean}
\end{table}

\section*{Declaration of competing interest}

The authors declare that they have no known competing financial interests or personal relationships that could have appeared to influence the work reported in this paper.

\newpage
\section*{References}
\nocite{*}
\bibliographystyle{elsarticle-harv}
\bibliography{test}

\begin{thebibliography}{51}
\expandafter\ifx\csname natexlab\endcsname\relax\def\natexlab#1{#1}\fi
\providecommand{\url}[1]{\texttt{#1}}
\providecommand{\href}[2]{#2}
\providecommand{\path}[1]{#1}
\providecommand{\DOIprefix}{doi:}
\providecommand{\ArXivprefix}{arXiv:}
\providecommand{\URLprefix}{URL: }
\providecommand{\Pubmedprefix}{pmid:}
\providecommand{\doi}[1]{\href{http://dx.doi.org/#1}{\path{#1}}}
\providecommand{\Pubmed}[1]{\href{pmid:#1}{\path{#1}}}
\providecommand{\bibinfo}[2]{#2}
\ifx\xfnm\relax \def\xfnm[#1]{\unskip,\space#1}\fi
\bibitem[{Akhtar and Mian(2018)}]{Akhtar2021Advances}
\bibinfo{author}{Akhtar, N.}, \bibinfo{author}{Mian, A.}, \bibinfo{year}{2018}.
\newblock \bibinfo{title}{Threat of adversarial attacks on deep learning in computer vision: A survey}.
\newblock \bibinfo{journal}{IEEE Access} \bibinfo{volume}{6}, \bibinfo{pages}{14410--14430}.
\bibitem[{Alexiadis et~al.(2004)Alexiadis, Colyar, Halkias, Hranac and McHale}]{NGSIMdata}
\bibinfo{author}{Alexiadis, V.}, \bibinfo{author}{Colyar, J.}, \bibinfo{author}{Halkias, J.}, \bibinfo{author}{Hranac, R.}, \bibinfo{author}{McHale, G.}, \bibinfo{year}{2004}.
\newblock \bibinfo{title}{The next generation simulation program}.
\newblock \bibinfo{journal}{Institute of Transportation Engineers. ITE Journal} \bibinfo{volume}{74}, \bibinfo{pages}{22}.
\bibitem[{Ban and Gruteser(2010)}]{Ban2010Mobile}
\bibinfo{author}{Ban, X.}, \bibinfo{author}{Gruteser, M.}, \bibinfo{year}{2010}.
\newblock \bibinfo{title}{Mobile sensors as traffic probes: addressing transportation modeling and privacy protection in an integrated framework}, in: \bibinfo{booktitle}{Traffic and Transportation Studies 2010}, p. \bibinfo{pages}{750–767}.
\bibitem[{Ban et~al.(2011)Ban, Hao and Sun}]{Ban2011Real}
\bibinfo{author}{Ban, X.J.}, \bibinfo{author}{Hao, P.}, \bibinfo{author}{Sun, Z.}, \bibinfo{year}{2011}.
\newblock \bibinfo{title}{Real time queue length estimation for signalized intersections using travel times from mobile sensors}.
\newblock \bibinfo{journal}{Transportation Research Part C: Emerging Technologies} \bibinfo{volume}{19}, \bibinfo{pages}{1133--1156}.
\bibitem[{Biggio et~al.(2013)Biggio, Nelson and Laskov}]{Biggio2013Poisoning}
\bibinfo{author}{Biggio, B.}, \bibinfo{author}{Nelson, B.}, \bibinfo{author}{Laskov, P.}, \bibinfo{year}{2013}.
\newblock \bibinfo{title}{Poisoning attacks against support vector machines}.
\newblock \bibinfo{journal}{arXiv:1206.6389 [cs, stat]} \URLprefix \url{http://arxiv.org/abs/1206.6389}. \bibinfo{note}{arXiv: 1206.6389}.
\bibitem[{Bungert et~al.(2021)Bungert, Raab, Roith, Schwinn and Tenbrinck}]{bungert2021clip}
\bibinfo{author}{Bungert, L.}, \bibinfo{author}{Raab, R.}, \bibinfo{author}{Roith, T.}, \bibinfo{author}{Schwinn, L.}, \bibinfo{author}{Tenbrinck, D.}, \bibinfo{year}{2021}.
\newblock \bibinfo{title}{Clip: Cheap lipschitz training of neural networks}, in: \bibinfo{booktitle}{Scale Space and Variational Methods in Computer Vision: 8th International Conference, SSVM 2021, Virtual Event, May 16--20, 2021, Proceedings}, \bibinfo{organization}{Springer}. pp. \bibinfo{pages}{307--319}.
\bibitem[{Burkard and Lagesse(2017)}]{Burkard2017Analysis}
\bibinfo{author}{Burkard, C.}, \bibinfo{author}{Lagesse, B.}, \bibinfo{year}{2017}.
\newblock \bibinfo{title}{Analysis of causative attacks against svms learning from data streams}, in: \bibinfo{booktitle}{Proceedings of the 3rd ACM on International Workshop on Security And Privacy Analytics}, pp. \bibinfo{pages}{31--36}.
\bibitem[{Chen et~al.(2017)Chen, Ban, Wang, Wang, Siddique, Fan and Lee}]{Chen2017Understanding}
\bibinfo{author}{Chen, C.}, \bibinfo{author}{Ban, X.J.}, \bibinfo{author}{Wang, F.}, \bibinfo{author}{Wang, J.}, \bibinfo{author}{Siddique, C.}, \bibinfo{author}{Fan, R.}, \bibinfo{author}{Lee, J.}, \bibinfo{year}{2017}.
\newblock \bibinfo{title}{Understanding GPS and Mobile Phone Data for Origin-Destination Analysis}.
\newblock \bibinfo{type}{Technical Report}.
\newblock \URLprefix \url{https://www.fhwa.dot.gov/planning/tmip/publications/other_reports/gps_mobile_data/index.cfm}. \bibinfo{note}{issue: FHWA-HEP-19-027}.
\bibitem[{Chen et~al.(2019)Chen, Xu, Wu, Yue, Yuan and Wang}]{Chen2019Deer}
\bibinfo{author}{Chen, J.}, \bibinfo{author}{Xu, H.}, \bibinfo{author}{Wu, J.}, \bibinfo{author}{Yue, R.}, \bibinfo{author}{Yuan, C.}, \bibinfo{author}{Wang, L.}, \bibinfo{year}{2019}.
\newblock \bibinfo{title}{Deer crossing road detection with roadside lidar sensor}.
\newblock \bibinfo{journal}{Ieee Access} \bibinfo{volume}{7}, \bibinfo{pages}{65944--65954}.
\newblock \bibinfo{note}{ISBN: 2169-3536 publisher: IEEE}.
\bibitem[{De~Cock(2011)}]{de2011ames}
\bibinfo{author}{De~Cock, D.}, \bibinfo{year}{2011}.
\newblock \bibinfo{title}{Ames, iowa: Alternative to the boston housing data as an end of semester regression project}.
\newblock \bibinfo{journal}{Journal of Statistics Education} \bibinfo{volume}{19}.
\bibitem[{Deng(2012)}]{deng2012mnist}
\bibinfo{author}{Deng, L.}, \bibinfo{year}{2012}.
\newblock \bibinfo{title}{The mnist database of handwritten digit images for machine learning research [best of the web]}.
\newblock \bibinfo{journal}{IEEE signal processing magazine} \bibinfo{volume}{29}, \bibinfo{pages}{141--142}.
\bibitem[{Dontchev and Rockafellar(2009)}]{Dontchev2009Implicit}
\bibinfo{author}{Dontchev, A.L.}, \bibinfo{author}{Rockafellar, R.T.}, \bibinfo{year}{2009}.
\newblock \bibinfo{title}{Implicit functions and solution mappings}. volume \bibinfo{volume}{543}.
\newblock \bibinfo{publisher}{Springer}.
\bibitem[{Facchinei and Pang(2003)}]{Facchinei2003Finite-dimensional}
\bibinfo{author}{Facchinei, F.}, \bibinfo{author}{Pang, J.S.}, \bibinfo{year}{2003}.
\newblock \bibinfo{title}{Finite-dimensional variational inequalities and complementarity problems}.
\newblock \bibinfo{publisher}{Springer}.
\bibitem[{Fiacco(1983)}]{Fiacco83}
\bibinfo{author}{Fiacco, A.V.}, \bibinfo{year}{1983}.
\newblock \bibinfo{title}{Introduction to sensitivity and stability analyses in nonlinear programming}.
\newblock \bibinfo{publisher}{Academics Press}.
\bibitem[{Gower(2018)}]{Gower2018Convergence}
\bibinfo{author}{Gower, R.M.}, \bibinfo{year}{2018}.
\newblock \bibinfo{title}{Convergence theorems for gradient descent}.
\newblock \bibinfo{journal}{Lecture notes for Statistical Optimization} \bibinfo{note}{Citation Key: gower2018convergence}.
\bibitem[{Guo and Ban(2024)}]{guo2024network}
\bibinfo{author}{Guo, Q.}, \bibinfo{author}{Ban, X.J.}, \bibinfo{year}{2024}.
\newblock \bibinfo{title}{Network multiscale urban traffic control with mixed traffic flow}.
\newblock \bibinfo{journal}{Transportation research part B: methodological} \bibinfo{volume}{185}, \bibinfo{pages}{102963}.
\bibitem[{Herrera et~al.(2010)Herrera, Work, Herring, Ban, Jacobson and Bayen}]{herrera2010evaluation}
\bibinfo{author}{Herrera, J.C.}, \bibinfo{author}{Work, D.B.}, \bibinfo{author}{Herring, R.}, \bibinfo{author}{Ban, X.J.}, \bibinfo{author}{Jacobson, Q.}, \bibinfo{author}{Bayen, A.M.}, \bibinfo{year}{2010}.
\newblock \bibinfo{title}{Evaluation of traffic data obtained via gps-enabled mobile phones: The mobile century field experiment}.
\newblock \bibinfo{journal}{Transportation Research Part C: Emerging Technologies} \bibinfo{volume}{18}, \bibinfo{pages}{568--583}.
\bibitem[{Huang et~al.(2011)Huang, Joseph, Nelson, Rubinstein and Tygar}]{Huang2011Adversarial}
\bibinfo{author}{Huang, L.}, \bibinfo{author}{Joseph, A.D.}, \bibinfo{author}{Nelson, B.}, \bibinfo{author}{Rubinstein, B.I.}, \bibinfo{author}{Tygar, J.D.}, \bibinfo{year}{2011}.
\newblock \bibinfo{title}{Adversarial machine learning}, in: \bibinfo{booktitle}{Proceedings of the 4th ACM workshop on Security and artificial intelligence}, pp. \bibinfo{pages}{43--58}.
\bibitem[{Iutzeler and Malick(2020)}]{Iutzeler2020Nonsmoothness}
\bibinfo{author}{Iutzeler, F.}, \bibinfo{author}{Malick, J.}, \bibinfo{year}{2020}.
\newblock \bibinfo{title}{Nonsmoothness in machine learning: specific structure, proximal identification, and applications}.
\newblock \bibinfo{journal}{Set-Valued and Variational Analysis} \bibinfo{volume}{28}, \bibinfo{pages}{661–678}.
\newblock \bibinfo{note}{Publisher: Springer}.
\bibitem[{Jagielski et~al.(2018)Jagielski, Oprea, Biggio, Liu, Nita-Rotaru and Li}]{Jagielski2018Manipulating}
\bibinfo{author}{Jagielski, M.}, \bibinfo{author}{Oprea, A.}, \bibinfo{author}{Biggio, B.}, \bibinfo{author}{Liu, C.}, \bibinfo{author}{Nita-Rotaru, C.}, \bibinfo{author}{Li, B.}, \bibinfo{year}{2018}.
\newblock \bibinfo{title}{Manipulating machine learning: Poisoning attacks and countermeasures for regression learning}, in: \bibinfo{booktitle}{2018 IEEE symposium on security and privacy (SP)}, \bibinfo{organization}{IEEE}. pp. \bibinfo{pages}{19--35}.
\bibitem[{Koh et~al.(2022)Koh, Steinhardt and Liang}]{Koh2022Stronger}
\bibinfo{author}{Koh, P.W.}, \bibinfo{author}{Steinhardt, J.}, \bibinfo{author}{Liang, P.}, \bibinfo{year}{2022}.
\newblock \bibinfo{title}{Stronger data poisoning attacks break data sanitization defenses}.
\newblock \bibinfo{journal}{Machine Learning} , \bibinfo{pages}{1–47}\bibinfo{note}{Publisher: Springer}.
\bibitem[{Li and Vorobeychik(2018)}]{Li2018Evasion-Robust}
\bibinfo{author}{Li, B.}, \bibinfo{author}{Vorobeychik, Y.}, \bibinfo{year}{2018}.
\newblock \bibinfo{title}{Evasion-robust classification on binary domains}.
\newblock \bibinfo{journal}{ACM Transactions on Knowledge Discovery from Data} \bibinfo{volume}{12}, \bibinfo{pages}{1--32}.
\newblock \bibinfo{note}{Number: 4}.
\bibitem[{Li et~al.(2016)Li, Wang, Singh and Vorobeychik}]{Li2016Data}
\bibinfo{author}{Li, B.}, \bibinfo{author}{Wang, Y.}, \bibinfo{author}{Singh, A.}, \bibinfo{author}{Vorobeychik, Y.}, \bibinfo{year}{2016}.
\newblock \bibinfo{title}{Data poisoning attacks on factorization-based collaborative filtering}.
\newblock \bibinfo{journal}{arXiv:1608.08182 [cs]} \URLprefix \url{http://arxiv.org/abs/1608.08182}. \bibinfo{note}{arXiv: 1608.08182}.
\bibitem[{Li et~al.(2021)Li, Liu, Dong, Zhao, Xue, Zhu and Lu}]{li2021hidden}
\bibinfo{author}{Li, S.}, \bibinfo{author}{Liu, H.}, \bibinfo{author}{Dong, T.}, \bibinfo{author}{Zhao, B.Z.H.}, \bibinfo{author}{Xue, M.}, \bibinfo{author}{Zhu, H.}, \bibinfo{author}{Lu, J.}, \bibinfo{year}{2021}.
\newblock \bibinfo{title}{Hidden backdoors in human-centric language models}, in: \bibinfo{booktitle}{Proceedings of the 2021 ACM SIGSAC Conference on Computer and Communications Security}, pp. \bibinfo{pages}{3123--3140}.
\bibitem[{Li and Ban(2019)}]{Li2019Connected}
\bibinfo{author}{Li, W.}, \bibinfo{author}{Ban, X.}, \bibinfo{year}{2019}.
\newblock \bibinfo{title}{Connected vehicles based traffic signal timing optimization}.
\newblock \bibinfo{journal}{IEEE Transactions on Intelligent Transportation Systems} \bibinfo{volume}{20}, \bibinfo{pages}{4354--4366}.
\newblock \bibinfo{note}{Event-title: IEEE Transactions on Intelligent Transportation Systems}.
\bibitem[{Li et~al.(2020)Li, Ban, Zheng, Liu, Gong and Li}]{li2020real}
\bibinfo{author}{Li, W.}, \bibinfo{author}{Ban, X.J.}, \bibinfo{author}{Zheng, J.}, \bibinfo{author}{Liu, H.X.}, \bibinfo{author}{Gong, C.}, \bibinfo{author}{Li, Y.}, \bibinfo{year}{2020}.
\newblock \bibinfo{title}{Real-time movement-based traffic volume prediction at signalized intersections}.
\newblock \bibinfo{journal}{Journal of Transportation Engineering, Part A: Systems} \bibinfo{volume}{146}, \bibinfo{pages}{04020081}.
\bibitem[{Liu et~al.(2022)Liu, Liu and Jiang}]{liu2022practical}
\bibinfo{author}{Liu, F.}, \bibinfo{author}{Liu, H.}, \bibinfo{author}{Jiang, W.}, \bibinfo{year}{2022}.
\newblock \bibinfo{title}{Practical {Adversarial} {Attacks} on {Spatiotemporal} {Traffic} {Forecasting} {Models}}.
\newblock \URLprefix \url{http://arxiv.org/abs/2210.02447}. \bibinfo{note}{arXiv:2210.02447 [cs]}.
\bibitem[{Liu et~al.(2023)Liu, Tian, Miranda-Moreno and Sun}]{liu2023adversarial}
\bibinfo{author}{Liu, F.}, \bibinfo{author}{Tian, J.}, \bibinfo{author}{Miranda-Moreno, L.}, \bibinfo{author}{Sun, L.}, \bibinfo{year}{2023}.
\newblock \bibinfo{title}{Adversarial {Danger} {Identification} on {Temporally} {Dynamic} {Graphs}}.
\newblock \bibinfo{journal}{IEEE Transactions on Neural Networks and Learning Systems} , \bibinfo{pages}{1--12}\URLprefix \url{https://ieeexplore.ieee.org/document/10068359/}.
\bibitem[{Liu et~al.(2011)Liu, Ning and Reiter}]{Liu2011False}
\bibinfo{author}{Liu, Y.}, \bibinfo{author}{Ning, P.}, \bibinfo{author}{Reiter, M.K.}, \bibinfo{year}{2011}.
\newblock \bibinfo{title}{False data injection attacks against state estimation in electric power grids}.
\newblock \bibinfo{journal}{ACM Trans. Inf. Syst. Secur.} \bibinfo{volume}{14}.
\newblock \URLprefix \url{https://doi.org/10.1145/1952982.1952995}. \bibinfo{note}{publisher-place: New York, NY, USA publisher: Association for Computing Machinery}.
\bibitem[{Luo et~al.(1996)Luo, Pang and Ralph}]{luo1996mathematical}
\bibinfo{author}{Luo, Z.Q.}, \bibinfo{author}{Pang, J.S.}, \bibinfo{author}{Ralph, D.}, \bibinfo{year}{1996}.
\newblock \bibinfo{title}{Mathematical programs with equilibrium constraints}.
\newblock \bibinfo{publisher}{Cambridge University Press}.
\bibitem[{Mei and Zhu(2015)}]{Mei2015Using}
\bibinfo{author}{Mei, S.}, \bibinfo{author}{Zhu, X.}, \bibinfo{year}{2015}.
\newblock \bibinfo{title}{Using machine teaching to identify optimal training-set attacks on machine learners}, in: \bibinfo{booktitle}{Proceedings of the aaai conference on artificial intelligence}.
\bibitem[{Mockus(1998)}]{mockus1998application}
\bibinfo{author}{Mockus, J.}, \bibinfo{year}{1998}.
\newblock \bibinfo{title}{The application of bayesian methods for seeking the extremum}.
\newblock \bibinfo{journal}{Towards global optimization} \bibinfo{volume}{2}, \bibinfo{pages}{117}.
\bibitem[{Papernot et~al.(2016)Papernot, McDaniel, Jha, Fredrikson, Celik and Swami}]{Papernot2016limitations}
\bibinfo{author}{Papernot, N.}, \bibinfo{author}{McDaniel, P.}, \bibinfo{author}{Jha, S.}, \bibinfo{author}{Fredrikson, M.}, \bibinfo{author}{Celik, Z.B.}, \bibinfo{author}{Swami, A.}, \bibinfo{year}{2016}.
\newblock \bibinfo{title}{The limitations of deep learning in adversarial settings}, in: \bibinfo{booktitle}{2016 IEEE European symposium on security and privacy (EuroS\&P)}, \bibinfo{organization}{IEEE}. pp. \bibinfo{pages}{372--387}.
\bibitem[{Perlin(1985)}]{perlin1985image}
\bibinfo{author}{Perlin, K.}, \bibinfo{year}{1985}.
\newblock \bibinfo{title}{An image synthesizer}.
\newblock \bibinfo{journal}{ACM Siggraph Computer Graphics} \bibinfo{volume}{19}, \bibinfo{pages}{287--296}.
\bibitem[{Polyak(1966)}]{polyak1966existence}
\bibinfo{author}{Polyak, B.T.}, \bibinfo{year}{1966}.
\newblock \bibinfo{title}{Existence theorems and convergence of minimizing sequences for extremal problems with constraints}, in: \bibinfo{booktitle}{Doklady Akademii Nauk}, \bibinfo{organization}{Russian Academy of Sciences}. pp. \bibinfo{pages}{287--290}.
\bibitem[{Reda et~al.(2021)Reda, Anwar, Mahmood and Tari}]{Reda2021Taxonomy}
\bibinfo{author}{Reda, H.T.}, \bibinfo{author}{Anwar, A.}, \bibinfo{author}{Mahmood, A.N.}, \bibinfo{author}{Tari, Z.}, \bibinfo{year}{2021}.
\newblock \bibinfo{title}{A taxonomy of cyber defence strategies against false data attacks in smart grid}.
\newblock \bibinfo{journal}{arXiv:2103.16085 [cs, eess]} \URLprefix \url{http://arxiv.org/abs/2103.16085}. \bibinfo{note}{arXiv: 2103.16085}.
\bibitem[{Suciu et~al.(2018)Suciu, Marginean, Kaya, Daume~III and Dumitras}]{suciu2018does}
\bibinfo{author}{Suciu, O.}, \bibinfo{author}{Marginean, R.}, \bibinfo{author}{Kaya, Y.}, \bibinfo{author}{Daume~III, H.}, \bibinfo{author}{Dumitras, T.}, \bibinfo{year}{2018}.
\newblock \bibinfo{title}{When does machine learning $\{$FAIL$\}$? generalized transferability for evasion and poisoning attacks}, in: \bibinfo{booktitle}{27th USENIX Security Symposium (USENIX Security 18)}, pp. \bibinfo{pages}{1299--1316}.
\bibitem[{Sun and Ban(2013)}]{Sun2013Vehicle}
\bibinfo{author}{Sun, Z.}, \bibinfo{author}{Ban, X.J.}, \bibinfo{year}{2013}.
\newblock \bibinfo{title}{Vehicle classification using gps data}.
\newblock \bibinfo{journal}{Transportation Research Part C: Emerging Technologies} \bibinfo{volume}{37}, \bibinfo{pages}{102--117}.
\bibitem[{Suya et~al.(2021)Suya, Mahloujifar, Suri, Evans and Tian}]{Suya2021Model-targeted}
\bibinfo{author}{Suya, F.}, \bibinfo{author}{Mahloujifar, S.}, \bibinfo{author}{Suri, A.}, \bibinfo{author}{Evans, D.}, \bibinfo{author}{Tian, Y.}, \bibinfo{year}{2021}.
\newblock \bibinfo{title}{Model-targeted poisoning attacks with provable convergence}, in: \bibinfo{booktitle}{International Conference on Machine Learning}, \bibinfo{organization}{PMLR}. pp. \bibinfo{pages}{10000--10010}.
\bibitem[{Tian et~al.(2022)Tian, Cui, Liang and Yu}]{tian2022comprehensive}
\bibinfo{author}{Tian, Z.}, \bibinfo{author}{Cui, L.}, \bibinfo{author}{Liang, J.}, \bibinfo{author}{Yu, S.}, \bibinfo{year}{2022}.
\newblock \bibinfo{title}{A comprehensive survey on poisoning attacks and countermeasures in machine learning}.
\newblock \bibinfo{journal}{ACM Computing Surveys} \bibinfo{volume}{55}, \bibinfo{pages}{1--35}.
\bibitem[{Tibshirani(1996)}]{Tibshirani1996Regression}
\bibinfo{author}{Tibshirani, R.}, \bibinfo{year}{1996}.
\newblock \bibinfo{title}{Regression shrinkage and selection via the lasso}.
\newblock \bibinfo{journal}{Journal of the Royal Statistical Society: Series B (Methodological)} \bibinfo{volume}{58}, \bibinfo{pages}{267–288}.
\newblock \bibinfo{note}{Publisher: Wiley Online Library}.
\bibitem[{Vorobeychik and Kantarcioglu(2018)}]{Vorobeychik2018Adversarial}
\bibinfo{author}{Vorobeychik, Y.}, \bibinfo{author}{Kantarcioglu, M.}, \bibinfo{year}{2018}.
\newblock \bibinfo{title}{Adversarial machine learning}.
\newblock \bibinfo{journal}{Synthesis Lectures on Artificial Intelligence and Machine Learning} \bibinfo{volume}{12}, \bibinfo{pages}{1--169}.
\newblock \bibinfo{note}{Publisher: Morgan \& Claypool Publishers}.
\bibitem[{Wang et~al.(2024a)Wang, Wang and Ban}]{wang2024datasurvey}
\bibinfo{author}{Wang, F.}, \bibinfo{author}{Wang, X.}, \bibinfo{author}{Ban, X.J.}, \bibinfo{year}{2024}a.
\newblock \bibinfo{title}{Data poisoning attacks in intelligent transportation systems: A survey}.
\newblock \bibinfo{journal}{Transportation Research Part C: Emerging Technologies} \bibinfo{volume}{165}, \bibinfo{pages}{104750}.
\bibitem[{Wang et~al.(2024b)Wang, Wang, Hong, Rockafellar and Ban}]{wang2024data}
\bibinfo{author}{Wang, F.}, \bibinfo{author}{Wang, X.}, \bibinfo{author}{Hong, Y.}, \bibinfo{author}{Rockafellar, R.T.}, \bibinfo{author}{Ban, X.J.}, \bibinfo{year}{2024}b.
\newblock \bibinfo{title}{Data poisoning attacks on traffic state estimation and prediction}.
\newblock \bibinfo{journal}{Transportation Research Part C: Emerging Technologies} , \bibinfo{pages}{104577}.
\bibitem[{Wang et~al.(2024c)Wang, Wang, Hong and Ban}]{wang2024transferability}
\bibinfo{author}{Wang, X.}, \bibinfo{author}{Wang, F.}, \bibinfo{author}{Hong, Y.}, \bibinfo{author}{Ban, X.}, \bibinfo{year}{2024}c.
\newblock \bibinfo{title}{Transferability in data poisoning attacks on spatiotemporal traffic forecasting models}.
\newblock \bibinfo{journal}{Available at SSRN 4827065} .
\bibitem[{Wang and Chaudhuri(2018)}]{Wang2018Data}
\bibinfo{author}{Wang, Y.}, \bibinfo{author}{Chaudhuri, K.}, \bibinfo{year}{2018}.
\newblock \bibinfo{title}{Data poisoning attacks against online learning} \URLprefix \url{http://arxiv.org/abs/1808.08994}. \bibinfo{note}{arXiv:1808.08994 [cs, stat]}.
\bibitem[{Xie et~al.(2022)Xie, Wang, Kong and Hong}]{Xie2022Universal}
\bibinfo{author}{Xie, S.}, \bibinfo{author}{Wang, H.}, \bibinfo{author}{Kong, Y.}, \bibinfo{author}{Hong, Y.}, \bibinfo{year}{2022}.
\newblock \bibinfo{title}{Universal 3-dimensional perturbations for black-box attacks on video recognition systems}.
\newblock \bibinfo{journal}{IEEE Security Privacy} , \bibinfo{pages}{19}.
\bibitem[{Xing et~al.(2013)Xing, Meng, Doozan, Snoeren, Feamster and Lee}]{xing2013take}
\bibinfo{author}{Xing, X.}, \bibinfo{author}{Meng, W.}, \bibinfo{author}{Doozan, D.}, \bibinfo{author}{Snoeren, A.C.}, \bibinfo{author}{Feamster, N.}, \bibinfo{author}{Lee, W.}, \bibinfo{year}{2013}.
\newblock \bibinfo{title}{Take this personally: Pollution attacks on personalized services}, in: \bibinfo{booktitle}{22nd USENIX Security Symposium (USENIX Security 13)}, pp. \bibinfo{pages}{671--686}.
\bibitem[{Yao et~al.(2024)Yao, Duan, Xu, Cai, Sun and Zhang}]{yao2024survey}
\bibinfo{author}{Yao, Y.}, \bibinfo{author}{Duan, J.}, \bibinfo{author}{Xu, K.}, \bibinfo{author}{Cai, Y.}, \bibinfo{author}{Sun, Z.}, \bibinfo{author}{Zhang, Y.}, \bibinfo{year}{2024}.
\newblock \bibinfo{title}{A survey on large language model (llm) security and privacy: The good, the bad, and the ugly}.
\newblock \bibinfo{journal}{High-Confidence Computing} , \bibinfo{pages}{100211}.
\bibitem[{Zhang and Zhang(2016)}]{zhang2016upper}
\bibinfo{author}{Zhang, Y.}, \bibinfo{author}{Zhang, L.}, \bibinfo{year}{2016}.
\newblock \bibinfo{title}{On the upper lipschitz property of the kkt mapping for nonlinear semidefinite optimization}.
\newblock \bibinfo{journal}{Operations Research Letters} \bibinfo{volume}{44}, \bibinfo{pages}{474--478}.
\bibitem[{Zhu et~al.(2021)Zhu, Feng, Pu and Ma}]{zhu2021adversarial}
\bibinfo{author}{Zhu, L.}, \bibinfo{author}{Feng, K.}, \bibinfo{author}{Pu, Z.}, \bibinfo{author}{Ma, W.}, \bibinfo{year}{2021}.
\newblock \bibinfo{title}{Adversarial {Diffusion} {Attacks} on {Graph}-based {Traffic} {Prediction} {Models}}.
\newblock \URLprefix \url{http://arxiv.org/abs/2104.09369}. \bibinfo{note}{arXiv:2104.09369 [cs]}.

\end{thebibliography}
\end{document}